\documentclass[11pt]{amsart}
\usepackage{amscd,amssymb}
\usepackage[mathscr]{eucal}

\title{dimensions of jet schemes of log singularities}
\author{Takehiko Yasuda }
\address{Department of Mathematical Sciences, University of Tokyo,
Komaba, Meguro, Tokyo, 153-8914, Japan }
\email{t-yasuda@ms.u-tokyo.ac.jp}

\theoremstyle{plain}
\newtheorem{thm}{Theorem}[section]

\newtheorem{lem}[thm]{Lemma}

\theoremstyle{definition}
\newtheorem{defn}[thm]{Definition}

\theoremstyle{remark}

\def\AA{\mathbb A}

\newcommand{\QQ}{\mathbb Q}
\newcommand{\ZZ}{\mathbb Z}
\newcommand{\LL}{\mathbb L}
\newcommand{\Zz}{\mathbb Z _{\geq 0}}

\newcommand{\cJ}{\mathscr J}

\newcommand{\cO}{\mathscr O}

\newcommand{\MM}{\mathbb M}
\newcommand{\Mhat}{\hat\mathbb M}

\newcommand{\Spec}{\mathrm{Spec}\,}

\newcommand{\Supp}{\mathrm{Supp}\,}

\newcommand{\Img}{\mathrm{Image}}

\newcommand{\codim}{\mathrm{codim}\,}

\newcommand{\sing}{\mathrm{sing}}

\newcommand{\red}{\mathrm{red}}
\newcommand{\reg}{\mathrm{reg}}
\newcommand{\const}{(\mathrm{const})}

\numberwithin{equation}{section}

\begin{document}

\maketitle

%

\begin{abstract}
We characterize Kawamata log terminal singularities and log canonical
singularities by dimensions of jet schemes. Our main result is 
Theorem \ref{thm-main}.
\end{abstract}

\section*{Introduction}

For birational geometry, it is natural to consider a pair $(X,D)$ of 
 a normal variety $X$ and a $\QQ$-divisor on it, and to consider 
singularities of pairs.
\textit{KLT (Kawamata log terminal)} and 
\textit{LC (log canonical)}
 singularities form important classes of log singularities.
They are defined by using a log-resolution and discrepancies.
When $X$ is $\QQ$-Gorenstein, we can naturally define
 KLT and LC also for a pair $(X,qY)$ where $Y$ is 
a closed subscheme of $X$
and $q \in \QQ_{>0}$.

In \cite{jet2}, Musta{\c{t}}{\u{a}} proposed a new point of view in the study
of log singularities. He characterized KLT and LC pairs $(X,qY)$ 
with $X$ smooth 
via dimensions of jet schemes of $Y$.
The aim of this short paper
 is to extend his result to the case $X$ is $\QQ$-Gorenstein and
to shorten his proof.
As he did, we also use the motivic integration as a main tool,
which is invented by Kontsevich 
\cite{Orsay} and 
extended to singular varieties by Denef and Loeser \cite{germs}.

\subsection*{Acknowledgements}
This paper forms a part of my master's thesis.
First, I would like to thank Yujiro Kawamata for his encouragement
and advise. 
I also thank Nobuyuki Kakimi for useful conversation.


\section{Motivic integration}

In this section, we review the motivic integration invented by Kontsevich
\cite{Orsay} and extended by Denef, Loeser \cite{germs}. 
Craw's paper \cite{craw} is a nice introduction. 

Let $k$ be an algebraically closed field of 
characteristic zero.

\subsection{Extended Grothendieck rings of varieties}

The \textit{Grothendieck ring of $k$-varieties}, denoted
$K_0(\mathrm{Var}/k)$, is
the abelian group generated by the isomorphism classes $[X]$ of 
$k$-varieties with the relations $[X]=[X\setminus Y]+[Y]$ if $Y$ is
a closed subvariety of $X$. The ring structure is defined by $[X][Y]=
[X\times Y]$. 
Let $\LL$ be the class $[\AA^1]$ of the affine line and let 
$\MM$ be the localization $K_0(\mathrm{Var}/k)[\LL^{-1}]$.
We define the \textit{dimension} of an element of $\MM$ by
\begin{align*}
   \dim ( \sum_i [X_i]\LL^{m_i} ) := \max_i ( \dim X_i + m_i ) ,
\end{align*}
with the convention $ \dim \emptyset := - \infty $.
Then, the following hold: for $\alpha,\ \beta \in \MM$,
\begin{itemize}
  \item $\dim \alpha \beta = \dim \alpha + \dim \beta$,
  \item $\dim (\alpha + \beta) \le \max \{ \dim \alpha , \dim \beta \} $ and 
     the equality holds if $\dim \alpha \ne \dim \beta$.
\end{itemize}
If we set $F_m \MM := \{ \alpha \in \MM  |  \dim \alpha \le m \}$,
then $(F_m \MM)_m$ 
is an ascending filtration with $F_m \MM \cdot F_n \MM \subset
F_{m+n} \MM$.
We define the \textit{complete Grothendieck ring of $k$-varieties} by
\[
\Mhat:= \lim_{\longleftarrow} \MM /F_m\MM \hspace{5mm} (m \to - \infty).
\]
Any element $\alpha$ of $\Mhat$ is expressed
as $\alpha = \sum _{m\in \ZZ} \alpha _m$ such that 
\begin{itemize}
	\item $\alpha _m = 0$ or $\dim \alpha _m = m$,
      \item $\alpha _m = 0$ for $m \gg 0$.
\end{itemize}
Then we define the \textit{dimension} of $\alpha$ to be 
$\max \{ m  |  \alpha _m \ne 0 \}$.

We denote by $ \Mhat ^ \QQ$ the ring $ \Mhat [\LL ^{q};q \in \QQ ] $.
Any element of $\Mhat ^ \QQ$ is expressed as $\alpha \LL^q$ with 
$\alpha \in \Mhat,\ q \in \QQ$.
We define the \textit{dimension}
of $\alpha \LL^q$ to be $ \dim \alpha + q $.

\subsection{Jet schemes and motivic measures}

Let $X$ be a $k$-scheme.  For $n \in \Zz \cup \{\infty\}$, 
an \textit{$n$-jet} on $X$ is a morphism
\[
   \Spec k[[t]]/(t^{n+1}) \to X ,
\]
where we used the convention $(t^{\infty})= (0) \subset k[[t]]$.
 The moduli scheme of the $n$-jets on $X$ always exists.
We call it the \textit{$n$-jet scheme} of $X$ and denote it by $L_n (X)$.
If $X$ is of finite type, 
then for $n < \infty$, so is $L_n (X)$.
If $X$ is smooth and of pure dimension $d$, then,
 for each $n \in \Zz$, the natural projection
$L_{n+1}(X)\to L_n(X)$ is a Zariski locally trivial $\AA^d$-bundle.

Now assume that $X$ is a $k$-variety of pure dimension $d$.
Let $\pi_n:L_\infty (X) \to L_n (X)$ be the canonical projection.

\begin{defn}
A subset $A$ of $L_\infty (X)$ is \textit{stable at level $n$} if we have:
\begin{enumerate}
	\item $\pi_n(A)$ is a constructible subset of $L_n (X)$,
        \item $A=\pi_n^{-1}\pi_n(A)$,
        \item for any $m\ge n$, the projection
$
\pi_{m+1}(A) \to \pi_m(A)
$
is a piecewise trivial $\AA^d$-bundle.
\end{enumerate}
(A morphism $f:Y\to X$ of schemes is called a \textit{piecewise trivial
$\AA^d$-bundle} if there is a stratification $X=\coprod X_i$ such that 
$f|_{f^{-1}(X_i)}:f^{-1}(X_i) \to X_i$ is isomorphic to 
$X_i\times \AA^d \to X_i$ for each i.)
A subset $A$ of $L_\infty (X)$ is \textit{stable} if it is 
stable at level $n$ for some $n\in \Zz$.
\end{defn}

We see that the stable subsets of $L_\infty (X)$ constitute a Boolean algebra,
that is, stable under finite intersection and finite union. 
For a stable
subset $A$ of $L_\infty (X)$,
an element $[\pi_n(A)]\LL^{-nd}\in\Mhat$ is independent of the choice of
$n \gg 0$. (We can define the class of a constructible set in 
$\MM$ and $\Mhat$ in the evident fashion.)
So the map
\begin{align*}
  \mu_X:\{\text{stable subsets of }L_\infty (X) \}&\to \Mhat \\
                      A\hspace{10mm}&\mapsto [\pi_n(A)]\LL^{-nd} \ (n \gg 0)
\end{align*}
is a finite additive measure. 
We can extend $\mu_X$ to the family of the \textit{measurable subsets} of 
$L_\infty X$, which is a family `big enough'. 
For details, see \cite{DL-quotient}, \cite{looi}.
We call $\mu_X$ the \textit{motivic measure} on $L_\infty X$. 

Let $A\subset L_\infty (X)$ be a measurable subset and
$\nu:A \to \QQ \cup \{ \infty \}$ a function.
We say that $\nu$ is a \textit{measurable
function} if 
the fibers 
are measurable and $\mu_X(\nu^{-1}(\infty)) = 0$. 

\begin{defn}
For a measurable function $\nu$, 
we formally define the
\textit{motivic integration} of $\LL^\nu$ by
\[
\int_A \LL^\nu d\mu_X := \sum_{n\in\QQ}\mu_X(\nu^{-1}(n))\LL^n .
\]
We say that $\LL^\nu$ is \textit{integrable} if this infinite sum converges 
in $\Mhat ^\QQ$.
\end{defn}

Let $Y \subset X$ be a closed subscheme and $\mathfrak a$ its ideal sheaf.
A closed point $\gamma \in L_\infty (X)$ corresponds to a morphism 
$\gamma' : \Spec k[[t]] \to X$. 
The function
\begin{align*}
F_Y :L_\infty (X) & \to \Zz\cup \{\infty\} \\
\gamma \; & \mapsto  n  \;  \text{ if }(\gamma')^{-1}\mathfrak a=(t^n) 
\end{align*}
is a measurable function.
For a $\QQ$-divisor $D=\sum_i q_i D_i$ on $X$ with $D_i$ a prime divisor, 
we define a measurable function $F_D := \sum _i q_i F_{D_i}$.
The following \textit{explicit formula} give a way to compute motivic
integrations.

\begin{lem}\label{lem-explicit}
Assume that $X$ is smooth and $\sum _{i=1}^s D_i$ is a SNC divisor on $X$ with
$D_i$ a prime divisor.

Let $ m_i \in \Zz, \ (1 \le i \le s) $,
and put $J := \{i| m_i >0 \} \subset \{1, \dots, s\}$.
Then we have
\[
   \mu_X ( \bigcap_i F_{D_i}^{-1}(m_i) ) = [D _J ^\circ] ( \LL - 1 )^{|J|} 
    \LL ^ {-\sum m_i} ,
\]
where
\[ 
D _J ^\circ := \bigcap _{i \in J} D _i  \bigg \backslash 
\bigcup _{ i \in \{  1, \dots , s  \} \setminus J}  D _i .
\]
\end{lem}

\begin{proof}
See \cite[the proof of Thm.\ 2.15]{craw}.
\end{proof}


\section{Main theorem}

Let $X$ be a normal variety of dimension $d$. 
Assume that $X$ is $\QQ$-Gorenstein, that is,
for some $r \in \ZZ_{>0} $,
$rK_X$ is a Cartier divisor. 
For a resolution $ p : \tilde X \to X$, the \textit{relative
 canonical divisor} $ K_{\tilde X /X} $  is
\[
   K_{\tilde X /X} := \frac{1}{r} ( r K_{\tilde X} - p^* (rK_X) ) .
\]

Let $Y$ be a closed subscheme of $X$ and $\mathfrak a$ its ideal sheaf. 
By Hironaka's theorem, there are a resolution $ p : \tilde X \to X $ and 
a SNC divisor $\sum _i  D_i$ on $\tilde X$  such
that 
\begin{itemize}
    \item $p^{-1}\mathfrak a = \cO_{\tilde X} (- \sum _i y_i D_i) $
     for some $y_i \in \Zz$,
    \item $K := K_{\tilde X /X} = \sum _i a_i D_i $ for some $a_i \in \QQ$.
\end{itemize}
Fix these notations through the rest of the paper.
In the proof of Theorem \ref{thm-main}, we suppose additional conditions
on $p$.

\begin{defn}
For $q \in \QQ_{>0}$, we say that the pair $ (X, qY) $ is 
\textit{KLT (Kawamata log terminal)},
resp.\ \textit{LC (log canonical)} if for every $i$, 
$ -q y_i + a_i + 1 > 0 $, resp.\ $ -q y_i + a_i + 1 \ge 0 $.
\end{defn}

Let $X_{\reg}$ be the smooth locus of $X$ and
 $\iota: X_{\reg} \hookrightarrow X$  the inclusion.
Then $ \omega^{[r]} _X = \iota_* ((\Omega _X^d)^{\otimes r})$
 is an invertible sheaf.
We define an ideal sheaf $\cJ \subset \cO_X$ by the following equation:
\[
   \cJ \omega^{[r]} _X = 
     \Img ( (\Omega _X ^d)^{\otimes r } \to \omega^{[r]} _X).
\]
Let $Z \subset X$  be the closed subscheme associated to $\cJ$.
Then $\Supp Z = X_{\sing}$.
The following is a variation of the transformation rule, a key of
the proof of the main theorem.

\begin{thm}\label{thm-trans}
Let $A \subset L_\infty X$ be a measurable subset and 
$\nu : A \to \QQ \cup \{ \infty \}$ a measurable function.
Then we have the following equality:
\[
   \int _ A \LL^{\nu + (1/r) \cdot F_Z } d\mu_X 
  = \int _{p_\infty^{-1}(A)} \LL^{\nu \circ p_\infty -F_K} d\mu_{\tilde X}.
\]
\end{thm}

\begin{proof}
It is a direct consequence of the transformation rule 
\cite[Lem.\ 3.3]{germs}.
\end{proof}

For each $e \in \Zz$, we put $A_e := F_Z^{-1} (e)$.
For each $n,e \in \Zz$, we define
\[
  L^e_n (Y) := L_n (Y) \cap \pi_n ( A_e ),
\]
where we take the intersection in $L_n (X)$.
Of course, $L^e_n (Y)$ depends on the inclusion $Y \subset X$.

\begin{lem}\label{lem-piecewise}
There is a positive integer $\theta$ such that for any $n,e\in \Zz$ with
$n \ge \theta e$, the natural projection 
$ \pi_{n+1} (A_e) \to \pi_{n} (A_e) $
is a piecewise trivial $\AA^d$-bundle.
\end{lem}

\begin{proof}
See \cite[Lem.\ 4.1]{germs}.
\end{proof}

\begin{thm}\label{thm-main}
Suppose $\Supp Z \supset X_{\sing}$. 
Let $\theta$ be a positive integer as in Lemma \ref{lem-piecewise}.
For $l \in \ZZ _{>0} $, 
let $lY$ be the closed subscheme of $X$ associated to
the ideal sheaf $\mathfrak a ^l$.
Suppose that $ \mathfrak a ^l \subset \cJ ^ \theta$.
Then,
\begin{enumerate}
 \item $ (X,qY) $ is KLT iff
      for any $e,n \in \Zz$ with $n \ge 
          \theta e$, 
    \begin{align*}
       \dim L^e_n(lY) + e/r  < (n +1 )( d-q/l  ). 
   \end{align*}
 \item $ (X,qY) $ is LC iff
      for any $e,n \in \Zz$ with $n \ge 
          \theta e$, 
    \begin{align*}
       \dim L^e_n(lY) + e/r  \le (n +1 )( d-q/l  ). 
   \end{align*}
 \end{enumerate}
\end{thm}

\begin{proof}
We prove only (1). 
For each $n \in \Zz$, let $B_n := F_Y^{-1}(n)$ and 
$B_{\ge n} := F_Y^{-1}( \ZZ _{\ge n}) $.
For each $e,n \in \Zz$ with $n \ge \theta e$, 
consider the following element of $\Mhat ^\QQ$:
\begin{align*}
 S(e,n)&:= \int _{B_{n+1} \cap A_e} \LL ^{qF_Y + (1/r) \cdot F_{Z} }
                   d\mu_X \\
    &= \mu_X ( B_{n+1} \cap A_e ) \LL ^{(n+1)q + e/r}.
\end{align*}
By Lemma \ref{lem-piecewise}, 
\begin{align*}
 S(e,n) = & (\mu_X (B_{\ge n + 1} \cap A_e ) 
   - \mu_X ( B_{ \ge n+2 })\cap A_e ))\LL ^{(n+1)q + e/r} \\
 = & ([L_n^e (Y)] - [L_{n+1}^e (Y)] \LL^{-d}) \LL ^{ -nd +(n+1)q + e/r } .
\end{align*}
By Lemma \ref{lem-piecewise} again, we have 
$\dim L_{n+1}^e (Y) \le \dim L_n^e (Y) + d $, and hence
\begin{equation}\label{ineq-1}
 \dim S(e,n) \le \dim L_n^e (Y) -nd +(n+1)q + e/r.
\end{equation}
If the equality does not hold, then 
$\dim L_{n+1}^e (Y) = \dim L_n^e (Y) + d $.

We suppose that 
$p^{-1}\cJ = \cO_{\tilde X}(- \sum _i z_iD_i)$ with $z_i > 0$. 

We set $\tilde Y:= \sum _i y_i D_i $ and $\tilde Z:= \sum _i z_i D_i $.
From Theorem \ref{thm-trans}, we have
\[
S(e,n) = \int \LL ^{ qF_{\tilde Y} -F_K } d\mu_{\tilde X}
\]
where the domain of the integration is
\[
 F_{\tilde Y}^{-1} (n+1) \cap F_{\tilde Z}^{-1}(e).
\]
From Lemma \ref{lem-explicit}, we obtain
\[
 S(e,n) = \sum _{J \subset \{ 1, \dots , s \}} \sum _{\mathbf m \in M} 
    [D_J^\circ] (\LL-1)^{|J|} \LL^{- \sum (-qy_i + a_i +1 )m_i},
\]
where 
\begin{align*}
M & = M (J , n , e )  \\
   & :=  \{ \mathbf m = (m_i)_{i\in J} \in (\ZZ_{>0})^J \ | \ 
          \sum y_i m_i = n \text{ and } \sum z_i m_i = e \}.
\end{align*}
For each $J$ and $\mathbf m$, the highest term of 
\[
   [D_J^\circ] (\LL-1)^{|J|} \LL^{- \sum (-qy_i + a_i +1 )m_i}
\]
has a positive coefficient. Hence,
\begin{equation}\label{eq-2}
 \dim S(e,n) = \max _{\substack{J \subset \{ 1,\dots, s \} 
\\ \mathbf m \in M}} ( d - \sum_i (-q y_i + a_i +1) m_i ) .
\end{equation}

`Only if' part: The proof is by contradiction.
So assume that $(X,qY)$ is KLT, that is, for every $i$, 
$-q y_i + a_i +1 > 0 $, and that for some $n',e' \in \Zz$ with
$n'>\theta e'$, \[
\dim L^{e'}_{n'} (Y) + e'/r \ge (n'+1)(d-q).
\]
By (\ref{eq-2}), we have $ \dim S(e',n') < d$.
On the other hand, since
\[
 \dim L^{e'}_{n'} (Y) -n'd +(n'+1)q + e'/r \ge d ,
\]
the equality in (\ref{ineq-1}) does not hold.
Hence $\dim L_{n'+1}^{e'} (Y) = \dim L_{n'}^{e'} (Y) + d $ and hence
$\dim L^{e'}_{n'+1} (Y) + e'/r \le ((n'+1)+1)(d-q)$. 
By the same argument, we have 
$\dim L_{n'+2}^{e'} (Y) = \dim L_{n'+1}^{e'} (Y) + d $
and so on. It contradicts Lemma \ref{lem-linear}.

`If' part:
From inequality (\ref{ineq-1}), we have that for any $e$, $n$ with 
$n \ge \theta e$,
\begin{equation}\label{ineq-3}
   \dim S(e,n) < d.
\end{equation}
In view of (\ref{eq-2}) and (\ref{ineq-3}), it is easy to find that 
for any $i$, $-q y_i + a_i +1 > 0$, that is, $(X,qY)$ is KLT.
We have thus completed the proof.
\end{proof}

\begin{lem}\label{lem-linear}
Then for each $e \in \Zz$, 
there is a real numbers $b$ such that $b < d$ and 
for every $n \in \Zz$,
\[
\dim  L_{n}^e (Y)  \le b  n+ \const .
\]
\end{lem}

\begin{proof}
It suffices to find an increasing
 linear function $\psi: \Zz \to \Zz$ and a real number
$b < d$ such that 
\[
 \dim L_{\psi (n)}^e (Y) \le b \psi (n) + \const .
\] 

First, the case $Y$  reduced:
For $m,n \in \Zz$ with $m \ge n$, we denote by $\pi_n^m $ the natural
projection $L_m(X) \to L_n(X)$.
By Greenberg's theorem \cite[Cor.\ 1]{greenberg3},
there is a linear function
$ g : \Zz \to \Zz $ such that 
\begin{itemize}
	\item for every $ n\in \Zz$, $g (n) \ge n$,
       \item  $ \pi_n ( L_\infty (Y) ) = \pi_n^{g(n)} (L_{g (n)} (Y))$.
\end{itemize}
 Since 
$ L^e_{g (n)} (Y) \subset (\pi_n^{g(n)})^{-1}
\pi_n^{g(n)} (L^e_{g(n)} (Y)) $,
\[
 \dim L^e_{g (n)} (Y) \le 
\dim (\pi_n^{g(n)})^{-1} \pi_n^{g(n)}(L^e_{g(n)} (Y)).
\]
By the definition of $g$, 
\[
 \pi_n^{g(n)} ( L_{g (n)}^e (Y) ) = \pi_n ( L^e_{\infty} (Y) ) .
\]
Let $d'$ be the dimension of $Y$.
By \cite[Lem 4.3]{germs}, 
\[
 \dim \pi_n^{g(n)} ( L_{g (n)}^e (Y) )=
\dim \pi_n ( L^e_{\infty} (Y) ) \le \dim (n+1)d' .
\]
Hence, by Lemma \ref{lem-piecewise}, 
\begin{align*}
 \dim L^e_{g (n)} (Y)  
 & \le (n+1)d' + (g (n) -n )d + \const\\
 & = g (n ) d + n ( d' -d ) + \const \\
 & = g ( n) ( d + (d'-d)/C_1) + \const ,
 \end{align*}
where $ C_1 $ is the constant such that 
$ g (n) = C_1 n + \const $. 
We have thus proved the assertion in this case.

The general case: It suffices to show the case
where $ Y = l (Y_{\red}) $ for some $l \in \ZZ _{>0} $.
By the definitions, we have the following
\begin{align*}
 L_{ln-1}^e (Y) & = \pi_{ln-1} ( F_{Y}^{-1} (ln) \cap A_e ) ,\\
 L_{n-1}^e (Y_{\red}) & = \pi _{n-1} ( F_{Y_{\red}}^{-1} (n) \cap A_e )  .
\end{align*}
Because for $n,e$ with $n \ge \theta e$, $F_{Y}^{-1} (ln) \cap A_e = F_{Y_{\red}}^{-1} (n) \cap A_e $
 is stable at level $n$, we have
 $L_{ln-1}^e (Y) = (\pi_{n-1}^{ln-1})^{-1} L_{n}^e (Y_{\red})$.
Therefore, for some $b<d$,
\begin{align*}
 \codim ( L_{l n-1}^e (Y) / \pi _ {l n-1} ( A_e ))  
 &=  \codim ( L_{n-1}^e (Y_{\red}) / \pi _ {n-1} ( A_e )) \\
 & \ge ( d -b ) \varphi (n)  + \const .
\end{align*}
Because $ \dim \pi_n (A_e) = dn + \const $, we have
\begin{align*}
\dim L_{l n -1}^e (Y) & \le d l n - ( d -b ) n + \const \\
 & \le (d - (d-b)/l) ( l n -1 ) + \const .
\end{align*}
This completes the proof.
\end{proof}


\end{document}